\documentclass[journal]{IEEEtran}
\IEEEoverridecommandlockouts
\usepackage{cite}
\usepackage{amsmath} 
\usepackage{amssymb}  
\usepackage{amsfonts}
\usepackage{mathrsfs}
\usepackage{mathtools}
\usepackage[dvipsnames]{xcolor}
\usepackage[thmmarks, amsmath]{ntheorem}
\usepackage{enumerate}
\usepackage{circuitikz}
\usepackage{pgfplots}
\usepackage{booktabs}
\usepackage[bbgreekl]{mathbbol}
\usepackage{colortbl}
\usepackage{ragged2e}
\usepackage{enumitem}
\usepackage{standalone}
\usepackage{tikz}
\usetikzlibrary{positioning, shapes}
\usetikzlibrary{decorations.pathreplacing}
\usetikzlibrary{arrows}

\usetikzlibrary{calc,trees,positioning,arrows,chains,shapes.geometric,%
    decorations.pathreplacing,decorations.pathmorphing,shapes,%
    matrix,shapes.symbols}

\tikzset{
>=stealth',
  punktchain/.style={rectangle, rounded corners, 
    draw=black, very thick,text width=10em, 
    minimum height=3em, text centered, on chain},
  line/.style={draw, thick, <-},
  element/.style={tape,top color=white,bottom color=blue!50!black!60!,
    minimum width=8em,draw=blue!40!black!90, very thick,
    text width=10em, minimum height=3.5em, text centered, on chain},
  every join/.style={->, thick,shorten >=1pt},
  tuborg/.style={decorate},
  tubnode/.style={midway, right=2pt},
}

\graphicspath{{./figures/}}
\DeclareMathOperator{\trace}{trace}

\DeclareSymbolFont{bbold}{U}{bbold}{m}{n}
\DeclareSymbolFontAlphabet{\mathbbold}{bbold}

\newcommand{\vzeros}[1][]{\mathbbold{0}_{#1}}

\newcommand{\tm}[1]{\textcolor{magenta}{#1}}
\newcommand{\tol}[1]{\textcolor{OliveGreen}{#1}}
\newcommand{\tor}[1]{\textcolor{Orange}{#1}}

\newcolumntype{R}[1]{>{\RaggedLeft\arraybackslash}p{#1}}

\newtheorem{rem}{Remark}

\usepackage{pgfplots}
\pgfplotsset{compat=1.6,
        scaled x ticks = false, 
        xticklabel style={/pgf/number format/fixed,/pgf/number format/precision=3},
        } 
\newlength\figureheight 
\newlength\figurewidth 

\renewcommand{\baselinestretch}{0.985}

\begin{document}
\tikzstyle{block} = [draw, fill=white, rectangle, minimum height=3em, minimum width=6em]
\tikzstyle{sum} = [draw, fill=white, circle, node distance=1cm]
\tikzstyle{input} = [coordinate]
\tikzstyle{output} = [coordinate]
\tikzstyle{pinstyle} = [pin edge={to-,thin,black}]
\tikzset{my cloud/.style={cloud, draw, aspect=2,cloud color={gray!5!white}}}
\tikzset{short/.append style={bipoles/length=.4cm, color = Black,very thick}}
\tikzset{european voltages}
\tikzset{american inductors}
\tikzset{american resistors}
\tikzset{american currents}

\title{Placement and Implementation of Grid-Forming and Grid-Following Virtual Inertia {and Fast Frequency Response}\thanks{This work was partially funded by the European Union's Horizon 2020 research and innovation programme under grant agreement N$^\circ$ 691800 and the SNF Assistant Professor Energy Grant \#160573. This article reflects only the authors' views and the European Commission is not responsible for any use that may be made of the information it contains.}}
\author{Bala~Kameshwar~Poolla,~\IEEEmembership{Student Member,~IEEE,}
	Dominic~Gro\ss{},~\IEEEmembership{Member,~IEEE,}
	and~Florian~D\"orfler,~\IEEEmembership{Member,~IEEE}
\thanks{B. K. Poolla, D.Gro\ss{}, and F. D{\"o}rfler are with the Automatic Control Laboratory, Swiss Federal Institute of Technology (ETH) Z\"urich, Switzerland. Email: {\tt \{bpoolla,grodo,dorfler\}@ethz.ch}.}}

\thispagestyle{plain}
\pagestyle{plain}
\maketitle

\begin{abstract}
The electric power system is witnessing a shift in the technology of generation. Conventional thermal generation based on synchronous machines is gradually being replaced by power electronics interfaced renewable generation. This new mode of generation, however, lacks the natural inertia and governor damping which are quintessential features of synchronous machines. The loss of these features results in increasing frequency excursions and, ultimately, system instability. Among the numerous studies on mitigating these undesirable effects, the main approach involves {\bf \emph {virtual}} inertia emulation to mimic the behavior of synchronous machines. In this work, explicit models of {\bf \emph{grid-following}} and {\bf \emph{grid-forming}} virtual inertia (VI) devices are developed for inertia emulation and fast frequency response in low-inertia systems. An optimization problem is formulated to optimize the parameters and location of these devices in a power system to increase its resilience. Finally, a case study based on a high-fidelity model of the South-East Australian system is used to illustrate the effectiveness of such devices. 
\end{abstract}

\begin{IEEEkeywords}
Low-inertia systems, Optimization methods, Power system dynamic stability.
\end{IEEEkeywords}
\IEEEpeerreviewmaketitle

\section{Introduction}

The past decade has seen a concerted focus on alternate sources of energy to replace conventional synchronous machine-based generation. A majority of the concerns forcing such a shift- namely greenhouse emissions, safety of nuclear generation and waste disposal, etc., are effectively addressed by cleaner alternatives, primarily-wind turbines and photovoltaics. These sources are interfaced by means of power electronic converters. Their large-scale integration, however, has raised concerns \cite{ENTSOE:16,AEMO:17,NordicReport:16} about system stability and especially frequency stability \cite{JGS-WLK:02, GL-JR-SR-DF-MJO:04,ASA-SR-GV-AC-DJH:17}. The inherent rotational inertia \cite{PT-DVH:16,WW-KE-GB-JK:15,FM-FD-GH-DH-GV:18} of the synchronous machines and the damping provided by governors assures system stability in the event of faults such as loss of generators, sudden fluctuation in power injections due to variable renewable sources, tie line faults, system splits, loss of loads, etc. In case of a frequency deviation, the inertia of synchronous machines acts as a first response by providing kinetic energy to the system (or absorbing energy). In contrast, converter interfaced generation fundamentally offers neither of these services, thus, making the system prone to instability.

Several studies have been carried out to propose control techniques to mitigate this loss of rotational inertia and damping. One extensively studied technique relates to using power electronic converters to mimic synchronous machine behavior \cite{QCZ-GW:11,CA-TJ-FD:17,AT-FD-FK-ZM-WH:18,HB-TI-YM:14}. These methods rely on concepts ranging from simple proportional-derivative to more complex controls under the name of, e.g., {\it Virtual Synchronous Generators}. All these strategies rely on some form of energy storage such as batteries, super-capacitors, flywheels, or the residual kinetic energy of wind turbines \cite{JM-SWHH-WLK-JAF:06}, which acts as a substitute for the kinetic energy of machines. 

These investigations have established the efficacy of virtual inertia (VI) and fast frequency response (FFR), i.e., primary frequency control without turbine delay, as a short-term replacement for machine inertia in low-inertia power systems. Also, as power converters operate at a much faster time scales compared to conventional generation, it is plausible to foresee future power systems based on predominantly converter-interfaced generation, without a major distinction between different time-scale controls such as inertia and fast frequency response, and primary frequency control provided by synchronous machines \cite{FM-FD-GH-DH-GV:18,JAT-SVD-DSC:16,DG-TP-MSD-FX-AM:17}. Here, we exclusively focus on power systems with reduced inertia due to loss of synchronous machines and utilize virtual inertia and fast primary frequency control as a remedy. 

Conventionally, the total inertia and primary frequency control in the system are the main metrics utilized for system resilience analysis \cite{NordicReport:16}. However, the authors in \cite{BKP-SB-FD:17} showed that not only is virtual inertia and primary frequency control vital, but its location in the power system is equally crucial and there can be a degradation in the performance due to ill-conceived spatial inertia distributions, even if the total virtual inertia added to the power system is identical \cite{DG-SB-BKP-FD:17}. Other commonly used performance metrics to quantify power system robustness include frequency nadir, {\it RoCoF} (Rate of change of frequency), and power system damping ratio \cite{AU-TSB-GA:14}. In the literature, the problem of optimally tuning and placing the virtual inertia and primary frequency controllers based on system norms (see \cite{MP-EH-BF-JWSP:17, BKP-SB-FD:17, DG-SB-BKP-FD:17, AM-UM-CH:16, FP-EM:17}) has been explored for small-scale test cases with linear models \cite{BKP-SB-FD:17, DG-SB-BKP-FD:17, AM-UM-CH:16}. In \cite{TSB-TL-DJH:15,TSB-FD:17,EDA-CZ-SG-SD-YC:18} time-domain and spectral metrics such as RoCoF, nadir, and damping ratios are considered. In \cite{TSB-TL-DJH:15, TSB-FD:17} a sequential linear programming approach is used to optimize the allocation of grid-following virtual inertia and primary frequency control. This method directly optimizes the frequency nadir and RoCoF. In \cite{EDA-CZ-SG-SD-YC:18} the power system is reduced to a single swing equation with first order turbine dynamics and the damping ratio and peak overshoot are optimized subject to an economic cost.

{As contributions, this paper develops explicit models of converter-based virtual inertia devices that capture the key dynamic characteristics of phase-locked loops (PLLs) \cite{C00} used in grid-following virtual inertia devices and of grid-forming controls such as virtual synchronous machines \cite{SD-JAS:13,QCZ-GW:11}, droop control \cite{FM-FD-GH-DH-GV:18,RO-UM-PA-GH:18}, and machine matching control \cite{AM14,CBB+15,CA-TJ-FD:17}, that can be used to provide virtual inertia. In addition, these models are suitable for integration with large-scale, non-linear power system models, thus allowing for parameter tuning through tractable optimization problems.

Moreover, the applicability of system norms as a performance metric for power system analysis is established beyond the prototypical swing equation, by considering detailed models. To this end, we propose a computationally efficient $\mathcal{H}_2$ norm based algorithm to optimally tune the parameters and the placement of the VI devices in order improve the resilience of low-inertia power systems. The key idea of this algorithm is to exploit the interpretation of VI devices as feedback controllers. Though the algorithm is applicable for a broader class of services offered by power electronic devices, we concentrate our analysis on virtual inertia and fast frequency response.

Finally, a high-fidelity model of the South-East Australian power system is modified to replicate a low-inertia scenario and used for an extensive case study. This modified system is augmented with VI devices to study their impact on system stability and validate the optimal tuning obtained by applying the proposed optimization algorithm. Moreover, through extensive simulations, we validate the linearized models used in the $\mathcal{H}_2$ optimization algorithm and study the impact of both grid-forming and grid-following virtual inertia on the disturbance responses of the non-linear power system. Lastly, time-domain simulations are presented to study system stability and to compare the response of grid-following and grid-forming virtual inertia in detail.
}

The remainder of the paper is structured as follows: In Section \ref{sec:model}, the power system model, converter models in both grid-following and grid-forming implementations are presented. The key performance metrics for grid stability and design constraints are identified and suitably defined in Section \ref{sec:perfconst}. In Section \ref{sec:metrics} a computational approach to identify the location of the inverters to improve post-fault response of the low-inertia power systems is proposed. In Section \ref{sec:test}, the low-inertia model based on the South-East Australian system is presented. A case study based on the two implementations of virtual inertia is presented in Section \ref{sec:results}, and suitable metrics are investigated to quantify the improvements in system stability. Finally, Section \ref{sec:conclusions} concludes the paper.

\section{System Model}
\label{sec:model}
We consider a high-fidelity, non-linear power system model, consisting of synchronous machines with governors, automatic voltage regulators (AVR), power system stabilizers (PSS), constant impedance and constant power loads, and renewable generation that is abstracted by constant power sources on the time scales of interest. The dynamical model of the power system is given by a differential-algebraic equation \cite{HM90}
\begin{subequations}
\label{eq:SysModel}
\begin{align}
\dot{x}_s&=f_s(x_s,z_s),\label{eq:SysModel:diff}\\
\vzeros[] &=g_s(x_s,z_s,{i}),\label{eq:SysModel:alg}
\end{align}
\end{subequations}
where $\vzeros[]$ is a vector of zeros, $x_s \in \mathbb{R}^{n_{x_s}}$ is the state vector and contains (but is not limited to) the mechanical states of the generators, their controllers and the states of other devices (i.e., non-linear dynamic loads and renewable in-feed). The three-phase transmission network is modeled by the algebraic equation \eqref{eq:SysModel:alg} in current-balance form \cite[Sec. 7.3.2]{PWS-MAP:97}. The vector $z_s \in \mathbb{R}^{n_{z_s}}$ comprises the AC signals of the transmission network, such as transmission line currents and bus voltages $v_k \in \mathbb{R}^3$. Moreover, the vector $i \in \mathbb{R}^{n_i}$ contains the three-phase currents $i_k \in \mathbb{R}^{3}$ injected at each bus $k \in \{1,\ldots,n_b\}$. In this section, we will use the input $i$ to incorporate explicit models of converter-based virtual inertia devices and disturbances into the power system model \eqref{eq:SysModel}. We note that, $z_s$ can be expressed as a function of the states $x_s$ and the current injections $i$. Moreover, after combining the power system model \eqref{eq:SysModel} with suitable models of power converters, we shall use the control inputs of the power converters to provide virtual inertia and fast frequency response. We now elaborate on the dynamics of the VI devices and the disturbance model.

\subsection{Modelling of virtual inertia devices}
The VI devices are power electronic devices that mimic the inertial response of synchronous generators. In the following we consider the two most common implementations- grid-following and grid-forming \cite{DG-TP-MSD-FX-AM:17}. A grid-following virtual inertia device is controlled to inject active power proportional to the frequency deviation and rate of change of frequency estimated by a phase-locked loop (PLL). In contrast, the grid-forming virtual inertia device is a voltage source that responds to power imbalances by changing the frequency of its voltage. In this paper, we model both types of VI devices as local dynamic feedback controllers. Even though we focus on two prototypical implementations of virtual inertia, the approach proposed in this paper can be used for any other arbitrary controller transfer function.
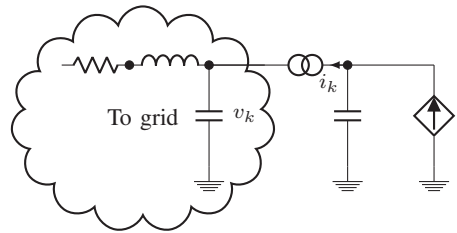
\begin{figure}[htp]
\centering
\begin{tikzpicture}[scale = 0.9,thick, every node/.style={font=\footnotesize, color = Black}]
\draw  [Black, thin] (0,0)   to [/tikz/circuitikz/bipoles/length=0.7cm, R, -]  ++(1,0) to[/tikz/circuitikz/bipoles/length=0.9cm,L=\(\), *-] ++(1.2,0) coordinate (p);
\draw [Black,thin, -] (p) -- (3,0);
\draw [Black,thin, -]  (5.5,-1.5) to[/tikz/circuitikz/bipoles/length=0.75cm, cI] (5.5, 0);
\draw [Black,thin, -] (4.2,0)--(5.5, 0);
\draw [Black,thin]   (p)++(-0.0234,0) to[/tikz/circuitikz/bipoles/length=0.6cm,C=$v_k$, *-] ++(0,-1.5);
\draw [Black,thin]   (p)++(2.0234,0) to[/tikz/circuitikz/bipoles/length=0.6cm,C, *-] ++(0,-1.5);
 \draw (2.1766,-1.5)  node [/tikz/circuitikz/bipoles/length=0.8cm, ground,thin] {};
\draw (4.2234,-1.5)  node [/tikz/circuitikz/bipoles/length=0.8cm, ground,thin] {};
\draw (5.5,-1.5)  node [/tikz/circuitikz/bipoles/length=0.8cm, ground,thin] {};
\draw [Black] (3.7,0) circle (0.15cm);
\draw [Black] (3.5,0) circle (0.15cm);
\draw [Black,thin, -] (p)-- (3.35,0);
\draw [Black,thin, -] (4.2,0) to (3.85,0);
\draw[Black,thin,-latex] (4.2,0) -- (3.95,0) node[below] {$i_k$};
\node [cloud, draw,cloud puffs=15,cloud ignores aspect, minimum width=3.5cm, minimum height=3cm] at (1.2,-0.8){\small{To grid}};
\end{tikzpicture}
\caption{Grid-following virtual inertia device.}
\label{Fig:GFoll-Conn}
\end{figure}

\paragraph {Grid-following} The grid-following VI device is controlled to inject active power proportional to the frequency deviation and RoCoF of the AC voltage at the bus where it is connected. To this end, each such virtual inertia device uses a frequency estimator, i.e., a Phase-Locked-Loop (PLL), that synchronizes to the bus voltage $v_k$ to obtain an estimate $\hat{\theta}_k$ of the bus voltage phase angle $\angle v_k$, the frequency estimate $\hat{\omega}_k$, and the RoCoF estimate $\dot{\hat{\omega}}_k$. We model such a synchronization device as
\begin{subequations}
\label{eq:VI_PLL}
\begin{align}
\dot{\hat{\theta}}_k &= \hat{\omega}_k,\\
\tau_k \dot{\hat{\omega}}_k &= - \hat{\omega}_k - K_\text{P,$k$} v_\text{q,$k$} - K_\text{I,$k$} \small{\int} v_\text{q,$k$},
\end{align}
\end{subequations}
where $v_\text{q,$k$}$ is the q-axis component of the bus voltage $v_k$ in a dq-frame with angle $\hat{\theta}_k$ and $v_\text{q,$k$} \approx \hat{\theta}_k - \angle v_k$ for small angle differences. Moreover, $\tau_k$, $K_\text{P,$k$}$, and  $K_\text{I,$k$}$ are the filter time constant, proportional gain, and integral synchronization gain. With $\tau_k=0$ the model \eqref{eq:VI_PLL} reduces to the standard synchronous reference frame phase locked loop (SRF-PLL) with a PI loop filter (i.e., $\dot{\hat{\theta}}_k = - K_\text{P,$k$} v_\text{q,$k$} - K_\text{I,$k$}  \int v_\text{q,$k$}$) commonly used in control of power converters \cite{GMF+14,C00}. However, the SRF-PLL with a PI loop filter does not provide an explicit RoCoF estimate. In contrast, incorporating a filter with time constant $\tau_k$ into the loop filter of the SRF-PLL allows us to obtain an explicit RoCoF estimate. 
\begin{rem}
Extensive simulations indicate that using a standard SRF-PLL in combination with a realizable differentiator results in worse control performance than integrating the RoCoF estimation into the PLL. Furthermore, the input into the PLL (i.e., $v_k$) is often subject to pre-filtering and $\tau_k$ can also be interpreted as moving the pre-filter into the loop filter (see the discussion in \cite{GMF+14}). 
\end{rem}
At the nominal steady-state, we have $\hat{\theta}_k \to \angle v_k$ and $\hat{\omega}_k \to 0$. This is because, we consider a reference frame rotating at the nominal grid frequency. With the frequency and RoCoF estimates in \eqref{eq:VI_PLL}, the VI device is modeled as
\begin{align}
P_{\text{VI},k}^\star& =K_{\text{\it foll}, k}\, [\hat{\omega}_k\,\,\dot{\hat{\omega}}_k]^\top, \quad Q_{\text{VI},k}^\star=0, \label{eq:VIFoll}
\end{align}
where $K_{\text{\it foll},k} =\begin{bmatrix}\tilde{d}_k\,\, \tilde{m}_k \end{bmatrix}$ are the control gains and $P_{\text{VI},k}^\star$ and $Q_{\text{VI},k}^\star$ are the set-points for the power injection of the grid-following VI device. The elements $\tilde{m}_k \geq 0$ are referred to as virtual inertia (reacts proportional to the derivative of the measured frequency), and $\tilde{d}_k \geq 0$ as the virtual damping (reacts proportional to the measured frequency itself). 

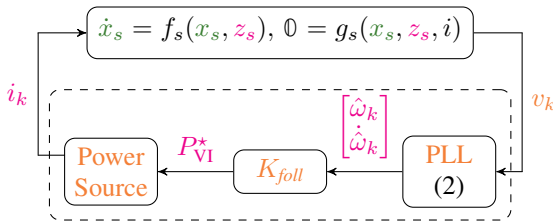
\begin{figure}[htp]
\centering
\begin{tikzpicture}
\tikzstyle{bigbox} = [draw, rounded corners, text width=4.9cm, align=center]
\tikzstyle{largebox} = [draw, rounded corners, text width=5.9cm, text height=1.5cm, align=center, dashed]
\tikzstyle{box} = [draw, rounded corners, text width=1cm, align=center]
\node[bigbox] (a) {\vspace{0.02in}$
\begingroup\tol{\dot{x}_s}\endgroup=f_s({ \tol {x_s}}, {\tm {z_s}}),\, \vzeros[] = g_s({ \tol {x_s}}, {\tm {z_s}}, i)$\vspace{0.02in}};
\node[left=0.65cm of a] (a2) {};
\node[right=0.65cm of a] (d1) {};
\node[box,below=1.2 cm of a] (b) {$\text{\tor{$K_\text{\it foll}$}}$};
\node[box,right=of b] (b1) {\tor {PLL} (2)};
\node[box,left=of b] (b2) {\tor{Power Source}};
\node[left=of b2] (b3) {};
\node[right=of b] (c) {};
\node[left=of b] (f) {};
\node[right=0.35cm of b1] (r) {};
\draw[->] (r) -- (b1);

\draw[->] (a2) -- (a);
\draw[-] (a) -- (d1);
\draw[-] (-3.217,0) -- (-3.217,-1.635) node[midway,left] {$\tm {i_k}$};
\draw[-] (-3.217,-1.628) -- (-2.86, -1.628);
\draw[->] (b) -- (f) node[midway, above] {$\tm {P_\text{VI}^\star}$};
\draw[->] (c) -- (b) node[midway, above=1pt] {$\tm{\begin{bmatrix} \hat\omega_{k}\\ \dot{\hat\omega}_{k} \end{bmatrix}}$};
\draw[-] (3.217, 0) -- (3.217,-1.85) node[midway,right] {$\tor{v_k} $};
\node[below=0.05cm of a] (b_l){};
\node[largebox, below=0.1 cm of b_l] (b_k){};
\end{tikzpicture}
\caption{Interconnection of a single grid-following virtual inertia device with power set-points according to \eqref{eq:VIFoll}.}
\label{Fig:Grid-Foll}
\end{figure}

The VI device utilizes a current source that injects the three-phase current $i_k$ at node $k$ (see Figure~\ref{Fig:GFoll-Conn}) and tracks the power references $P_{\text{VI},k}^\star$ and $Q_{\text{VI},k}^\star$ with time constant $\tau_\text{\it foll}=100$ ms. Figure~\ref{Fig:Grid-Foll} shows the overall control strategy. 

\begin{figure}[htp]
\centering
\begin{tikzpicture}[scale = 0.9,thick, every node/.style={font=\footnotesize, color =Black}]
\draw  [Black, thin] (0,0)   to [/tikz/circuitikz/bipoles/length=0.7cm, R, -]  ++(1,0) to[/tikz/circuitikz/bipoles/length=0.9cm,L=\(\), *-] ++(1.2,0) coordinate (p);
\draw [Black,thin, -] (p) -- (2,0);
\draw [Black,thin]   (p) to[/tikz/circuitikz/bipoles/length=0.7cm,sV<=$v_{\mathrm{VI},k}$] ++(0,-1.5);
\draw [Black,thin]   (p)++(-2.2234,0) to[/tikz/circuitikz/bipoles/length=0.7cm,C, *-] ++(0,-1.5);
\draw (2.2,-1.5)  node [/tikz/circuitikz/bipoles/length=0.8cm, ground,thin] {};
 \draw (-.0234,-1.5)  node [/tikz/circuitikz/bipoles/length=0.8cm, ground,thin] {};
\draw [Black,thin, -] (0,0) to (-0.55,0);
\draw[Black,thin,-latex] (0,0) -- (-0.35,0) node[below] {$i_k$};
\draw [Black,thin, -] (-2.5,0)-- (-1.05,0);
\draw [Black] (-0.7,0) circle (0.15cm);
\draw [Black] (-0.9,0) circle (0.15cm);
\draw  [Black, thin] (-4.5,0)   to [/tikz/circuitikz/bipoles/length=0.7cm, R, -]  ++(1,0) to[/tikz/circuitikz/bipoles/length=0.9cm,L=\(\), *-] ++(1.2,0) coordinate (q);
\draw [Black,thin]   (q)++(0.0234,0) to[/tikz/circuitikz/bipoles/length=0.7cm,C=$v_k$, *-] ++(0,-1.5);
\draw (-2.2766,-1.5)  node [/tikz/circuitikz/bipoles/length=0.8cm, ground,thin] {};
\node [cloud, draw,cloud puffs=15,cloud ignores aspect, minimum width=3.5cm, minimum height=3cm] at (-3.2,-0.8){\small{To grid}};
\end{tikzpicture}
\caption{Grid-forming virtual inertia device.}
\label{Fig:GForm-Conn}
\end{figure}
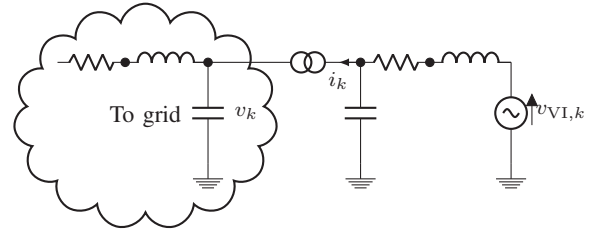

\begin{figure}[htp]
\centering
\begin{tikzpicture}
\tikzstyle{bigbox} = [draw, rounded corners, text width=4.9cm, align=center]
\tikzstyle{largebox} = [draw, rounded corners, text width=5cm, text height=1.1cm, align=center, dashed]
\tikzstyle{box} = [draw, rounded corners, text width=1cm, align=center]
\node[bigbox] (a) {\vspace{0.02in}$
\begingroup\tol{\dot{x}_s}\endgroup=f_s({ \tol {x_s}}, {\tm {z_s}}),\, \vzeros[] = g_s({ \tol {x_s}}, {\tm {z_s}}, i)$\vspace{0.02in}};
\node[left=0.3cm of a] (a2) {};
\node[right=0.3cm of a] (d) {};
\node[right=of a.350] (d1) {};
\node[below=1.2cm of a] (b){};
\node[below=1.164cm of a] (e){};
\node[box,left=of b] (q) {\tor{RLC filter}};
\draw[-] (a) to (e);
\node[right=of q.195] (l) {};
\node[right=of q.165] (m) {};
\node[right=1.12 cm of m] (v) {};
\draw[<-] (m)--(v) node[midway,above] {$v_k$};
\node[box,right=2.25cm of l] (b1) {$\tor{\text{VI}}\, (4)$};
\node[right=0.5cm of b1] (rb) {};
\draw[->]  (rb) -- (b1);
\node[left=of q] (b3) {};
\node[right=of b] (c) {};
\node[left=0.5cm of q] (f) {};
\draw[-] (q) -- (f);
\draw[->] (a2) -- (a);
\draw[-] (-2.87,0) -- (-2.87,-1.68) node[midway,left] {$\tm{i_k}$};
\draw[->] (b1) -- (l) node[midway,below] {$v_{\text{VI}, k}$};
\draw[-] (2.867,0) -- (2.867,-1.8495) node[midway,right] {$\tor{P_{\text{VI},k}} $};
\draw[-] (d) -- (a) ;
\node[below=0.35cm of a] (b_l){};
\node[largebox, below=0.1 cm of b_l] (b_k){};
\end{tikzpicture}
\caption{Interconnection of a single grid-forming VI device modeled via \eqref{eq:VIForm}}
\label{Fig:Grid-Form}
\end{figure}
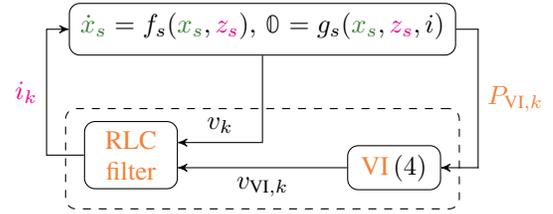

\paragraph {Grid-forming} The grid-forming VI device uses a voltage source connected to the grid via an LC filter with parasitic losses (see Figure~\ref{Fig:GForm-Conn}) 
%
%
that generates a voltage $v_{\text{VI},k}$ with an angle $\theta_{\text{VI},k} = \angle v_{\text{VI},k}$ that is a function of the power in-feed of the VI devices. The device is modeled via 
\begin{subequations}
\label{eq:VIForm}
\begin{align}
\dot{\theta}_{\text{VI},k} &= \omega_{\text{VI}, k},\\
\tilde{m}_k \dot{\omega}_{\text{VI},k} &= - \tilde{d}_k \omega_{\text{VI},k} - P_{\text{VI}, k}, \label{eq:VIForm:freq}
\end{align}
\end{subequations}
where $\theta_{\text{VI},k}$, $\omega_{\text{VI},k}$ are the angle and frequency of voltage generated by the grid-forming VI device, $P_{\text{VI},k}$ is the active power from the grid-forming VI device into the grid, $\tilde{m}_k > 0$ is the virtual inertia constant, and $\tilde{d}_k \geq 0$ the virtual damping constant. The amplitude of the voltage is regulated at the nominal operating voltage of the bus to which the device is connected. The overall signal flow for the grid-forming VI device is shown in Figure~\ref{Fig:Grid-Form}. 
%
%
The second-order active power-frequency droop characteristics \eqref{eq:VIForm} are the core operating principle of a wide range of grid-forming control algorithms for power converters. For instance, under the assumption that the controlled internal dynamics of the converter are sufficiently fast, droop control with a low pass filter in the power controller (see Figure 5 of \cite{PPG07}) is equivalent to \eqref{eq:VIForm} (see Lemma 4.1 of \cite{SG+13}). Similarly, \eqref{eq:VIForm} can be explicitly recovered for a wide range of grid-forming controls by applying model reduction techniques that eliminate the fast controlled internal dynamics of power converters \cite{CGD17}. This includes virtual synchronous machines that directly enforce a second-order frequency droop behavior (see Section II-A of \cite{SD-JAS:13} and Section II-B of \cite{QCZ-GW:11}) as well as controls based on matching the dynamics of power converters to that of a synchronous machine \cite{AM14,CBB+15,CA-TJ-FD:17,CGD17} (see Remark 2 and Section 3.3 of \cite{CA-TJ-FD:17}). 

\subsection{Disturbance model}\label{sec:dist}
We consider a general class of disturbance signals $\eta_k(t)$ that act at the voltage buses of the power system \eqref{eq:SysModel} through the current injection $i_k$. This approach can be used to model a wide range of faults such as load steps, fluctuations in renewable generation, or generator outages (i.e., by canceling the current injection of a generator). For brevity of presentation, we focus on faults that map changes in active power injection at every bus (i.e., changes in demand or generation) to current injections $i_k$ at every voltage bus $k$. We denote by $\eta = (\eta_1,\ldots,\eta_{n_d})$ the disturbance vector that corresponds to, e.g., changes in load or generation such as fluctuations of renewables.

\section{Performance metrics and design constraints}\label{sec:perfconst}
In this section, we discuss several performance metrics typically utilised in stability analysis of power systems. As an alternative to these conventional metrics, we propose system norms as a tool to assess transient stability and for optimizing the allocation of VI devices. Further, we discuss design constraints on the virtual inertia and damping gains arising from limits on the maximum power output of the VI devices and specifications imposed by grid-codes.

\subsection{Performance metrics}
\label{subsec:PerfMets}
Based on the model presented in Section~\ref{sec:model}, we now formally define a set of performance metrics that we shall use to assess the frequency stability of the grid, when subjected to disturbances.
Using the response of the system following a disturbance input $\eta (t)$ several time-domain metrics can be defined. In particular, given a negative step disturbance, e.g., a sudden load increase or loss of generation, at time $t=0$, we define the following indices on the time-domain evolution of the frequency vector $\omega = (\omega_1,\ldots,\omega_n)$ that collects the frequencies at different buses. Let the frequency nadir $\left\vert{\omega}_k\right\vert_\infty$ and maximum RoCoF\footnote{The RoCoF for the generators is computed by filtering the frequency derivative through a low-pass filter \cite{ENTSOE:15}.} $\left\vert\dot\omega_k\right\vert_\infty$ at each bus $k$ be given by
\begin{align}
\left\vert\dot\omega_k\right\vert_{\infty} &\coloneqq \max_{t\ge 0} \left\vert \dot \omega_k(t) \right\vert,
\label{eq: def rocof}\\
\left\vert{\omega_k} \right\vert_\infty &\coloneqq \left\vert \min_{t\ge 0}\, \omega_k(t) \right\vert.
\label{eq: def nadir}
\end{align}
Further, let $\omega_\text{G}=\begin{pmatrix}\omega_{\text{G},1}, \omega_{\text{G},2}, \ldots\end{pmatrix} = [\omega^\top_{\text{G},1},\omega^\top_{\text{G},2},\ldots]^\top$, $\dot\omega_\text{G}=\begin{pmatrix}\dot\omega_{\text{G},1},\dot\omega_{\text{G},2},\ldots \end{pmatrix}$, $P_\text{G}=\begin{pmatrix} P_{\text{G},1},P_{\text{G},2},\ldots\end{pmatrix}$, and $P_\text{VI}=\begin{pmatrix} P_{\text{VI},1},P_{\text{VI},2},\ldots\end{pmatrix}$ collect the generator frequencies, the RoCoF, the mechanical power injections, and the active power injections from VI devices. For the same step disturbance as considered above, we also define the peak power injection by the virtual inertia devices as well as the peak power injection due to the governor response of synchronous machines
\begin{align}
\!\left\vert P_{\text{VI},k} \right\vert_{\infty}\! \coloneqq \max_{t\ge 0} \left|P_{\text{VI},k}(t)\right|, \quad 
\left\vert P_{\text{G},k} \right\vert_{\infty}\! \coloneqq \max_{t\ge 0} \left|P_{\text{G},k}(t)\right|.\!
\label{eq: def peak}
\end{align}
Next, for any disturbance input $\eta$ (not necessarily a step) we define the metrics quantifying the energy imbalance and control effort on a time horizon $\tau$. The integrals
\begin{align}
 E_{\tau, \eta(t), \omega} \coloneqq \int_0^\tau \!\! \omega_\text{G}^\top \omega_\text{G}\,\text{d}t,\quad
 E_{\tau, \eta(t), \dot\omega} \coloneqq \int_0^\tau \!\! \dot{\omega}_{\mathrm{G}}^\top \,\dot{\omega}_\mathrm{G}\,\text{d}t,
\label{eq: def h2freq}
\end{align}
capture the frequency and RoCoF imbalance post-fault, and
\begin{align}
\!\!E_{\tau, \eta(t), P_\mathrm{VI}} \coloneqq \!\int_0^\tau \!\! P_\text{VI}^\top P_\text{VI} \,\text{d}t, \quad%
E_{\tau, \eta(t), P_\text{G}} \coloneqq \!\int_0^\tau \!\! P_\text{G}^\top P_\text{G} \,\text{d}t,\!
\label{eq: def h2effort}
\end{align}
encode the virtual inertia and damping effort of the converters, and the generator mechanical efforts. 

Consider a weighted sum of the metrics in \eqref{eq: def h2freq}-\eqref{eq: def h2effort}, i.e.,
\begin{align*}
J_{\tau, \eta(t)}\!\coloneqq\!\int_0^{\tau} \!\!\! r_\text{G} P_\text{G}^\top P_\text{G}\!+r_\text{VI} P_\text{VI}^\top P_\text{VI}+r_{\omega}  \omega_\text{G}^\top \omega_\text{G} +r_{\dot\omega}\, \dot{\omega}_{\mathrm{G}}^\top \dot{\omega}_\mathrm{G} \,\mathrm{d}t,
\end{align*}
with non-negative scalars $r_\text{G}$, $r_\omega$, $r_{\dot\omega}$, and $r_\text{VI}$ trading off the relative efforts. 

In the context of transient frequency stability analysis of power systems, faults are typically modeled as steps capturing, e.g., an increase in load or loss of generation. The main purpose of virtual inertia and fast frequency response is to improve the transient behavior of the power system immediately after such a fault. However, unless the time horizon $\tau$ is chosen carefully, the quadratic cost $J_{\tau, \eta(t)}$ for a step disturbance is a questionable metric for optimizing virtual inertia and fast frequency response because it is dominated by the post-fault steady-state deviation. To avoid this problem, we consider a metric referred to as the $\mathcal H_2$ norm that can be interpreted as the energy of the system response to impulse disturbance inputs. Specifically, the $\mathcal{H}_2$ norm can be obtained by perturbing the system with a unit impulse $\delta(t)$ at every disturbance input $\eta_k$ individually and summing over the resulting infinite horizon costs $J_{\infty, \delta_k(t)}$. To clarify the interpretation in the context of power systems, note that the impulse response of a linear system is equal to the time derivative of its step response. In other words, as time goes to infinity the step response of a stable system will settle to constant values, but the signals in the integral of the cost $J_{\infty, \delta_k(t)}$ tend to zero. Therefore, the $\mathcal{H}_2$ norm predominantly captures the initial transient. This is in line with the use of virtual inertia and fast frequency response to stabilize the frequency before slower controls and ancillary services act to control the long-term post-fault steady-state behavior of the system.

As the $\mathcal{H}_2$ norm measures short-term energy imbalance, it is a suitable proxy for transient power system stability. Concurrently, $\mathcal{H}_2$ norms result in tractable, well understood design and optimization problems that apply to a broader class of disturbances than the classical power system metrics and, in special cases, allow to solve min-max problems arising by maximizing over the disturbance vector $\eta$ while minimizing over inertia coefficients \cite{BKP-SB-FD:17}. In the remainder, we will use $\mathcal H_2$ norms for control design and tuning. However, for evaluation purposes, we also consider the metrics \eqref{eq: def rocof}--\eqref{eq: def peak} commonly used in power system analysis.

\subsection{Design constraints}\label{sec:constraint}
In addition to the performance metrics presented in the previous section, the virtual inertia devices are also subject to constraints on their power injection. Moreover, the net damping is constrained by grid-codes and primary control reserve markets. To account for constraints on the net damping, we impose an upper bound on the sum of damping gains of the VI devices, i.e., $\sum\nolimits_k \tilde{d}_k \leq d_{\text{sum}}$.
This ensures realistic results in line with power system operation. Further, we use constraints on the individual inertia and damping gains to account for the maximum power rating of the power converters. Notably, it has been observed from empirical data, that the maximum frequency deviation and the maximum RoCoF do not occur at the same time (see Section III-C of \cite{DG-SB-BKP-FD:17}, see also the scatter plot in Figure 3 of \cite{BKP-DG-TB-SB-FD:18}). Therefore, the inertia response and the damping response do not attain their peak values simultaneously. Based on this observation, the maximum power injection constraint can be approximated by the constraints 
\[ \tilde{m}_k \leq \dfrac{P_{\max,k}}{|\dot{\omega}|_{\max}}, \quad \tilde{d}_k \leq \dfrac{P_{\max,k}}{|\omega|_{\max}},\] 
where $P_{\max,k}$ is the power rating of the $k$-th converter and ${|\dot{\omega}|}_{\max}$ and ${|\omega|}_{\max}$ are a priori estimates of the maximum RoCoF and frequency deviation. In addition, we limit the individual damping and inertia gains to be non-negative.

\section{Closed-loop system and $\mathcal{H}_2$ optimization}
\label{sec:metrics}
In this section we present a computational approach to answer the question of ``how and where to optimally use virtual inertia and damping?'' via appropriate tuning of the gain matrices $K_\text{\it foll}$ and $K_\text{\it form}$ in order to improve the post-fault response of a low-inertia power system.

\subsection{Closed-loop system model and linearization}\label{sec:clandlin}
The placement and tuning of VI devices can be recast as a system norm (input-output gain) minimization problem for a linear system. To this end, we combine the power system model \eqref{eq:SysModel} with the disturbance model presented in Section \ref{sec:dist} and either grid-following \eqref{eq:VIFoll} or grid-forming \eqref{eq:VIForm} virtual inertia device models. Next, we define inputs and outputs of the system interconnected with the grid-forming and grid-following devices that allow us to optimize the virtual inertia and damping gains, $\tilde{m}_k$ and $\tilde{d}_k$ respectively. 

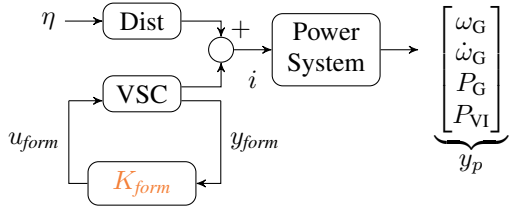
\begin{figure}[htp]
\centering
\begin{tikzpicture}
\tikzstyle{bigbox} = [draw, rounded corners, text width=0.8cm, align=center]
\tikzstyle{bigbox1} = [draw, rounded corners, text width=1.2cm, align=center]
\tikzstyle{box} = [draw, rounded corners, text width=1cm, align=center]
\tikzstyle{largebox} = [draw, rounded corners, text width=4.4cm, text height=0.7cm, align=center, dashed]
\node[bigbox1] (a) {Power\\ System};

\node[sum, left=0.5cm of a] (sum) {};

\node[left=1cm of a.205] (al){};
\node[right=0.5cm of a.345] (ar){};

\node[left=1cm of a.155] (kl){};
\node[right= of kl](rkl){};
\node[left=0.5cm of ar] (abr){};
\node[right=1cm of al] (abl){};
\node[below=0.00cm of abl] (abll){};
\node[bigbox, left=1.2 cm of abll] (al_VSC) {VSC};
\node[bigbox1,below=0.65cm of al_VSC] (b) {\tor{$K_\text{\it form}$}};

\node[left=1.8cm of b] (bl){};
\node[right=1.8cm of b] (br){};

\node[bigbox, above=0.475 cm of al_VSC] (al_d) {Dist};
\node[left=0.5cm of al_d](lkl){$\eta$};
\node[below=0.2cm of lkl](lklt){};
\node[left=0.2cm of al_VSC.190](lklb){};
\node[right=0.275cm of al_VSC.350](lklbr){};

\draw[->] (sum)-- (a)  node[midway,below=0.15 cm]{$i$};

\draw[->] (al_VSC.370)-|(sum);
\draw[->] (al_d)-|node[near end, right] {$+$} (sum);

\draw[-] (b.west)-|(lklb.west) node[near end,left]{$u_\text{\it form}$};
\draw[->] (lklb.west)|-(al_VSC.190);

\draw[-] (al_VSC.350)-|(lklbr.east);
\draw[->] (lklbr.east) |-(b.east)  node[near start,right]{$y_\text{\it form}$};

\draw[->] (lkl)-- (al_d);

\node[right=0.5cm of a](p){};
\node[below=0.3cm of p](x){};
\node[right=-0.15cm of x](y){$\underbrace{\begin{bmatrix}\omega_\text{G}\\ \dot{\omega}_\text{G}\\P_\text{G}\\ P_\text{VI} \end{bmatrix}}_{\displaystyle{y_p}}$};
\draw[->] (a)-- (p);

\end{tikzpicture}
\caption{Closed-loop system interconnection for the grid-following VI with tuning parameters $K_\text{\it form}$, where VSC is the voltage source converter.}
\label{Fig:ClosedSys-Form}
\end{figure}

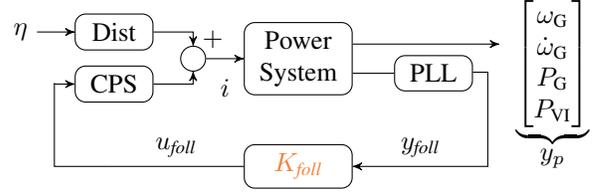
\begin{figure}[htp]
\centering
\begin{tikzpicture}
\tikzstyle{bigbox} = [draw, rounded corners, text width=0.8cm, align=center]
\tikzstyle{bigbox1} = [draw, rounded corners, text width=1.2cm, align=center]
\tikzstyle{box} = [draw, rounded corners, text width=1cm, align=center]
\tikzstyle{largebox} = [draw, rounded corners, text width=4.4cm, text height=0.7cm, align=center, dashed]
\node[bigbox1] (a) {Power\\ System};
\node[bigbox1,below=0.65cm of a] (b) {\tor{$K_\text{\it foll}$}};
\node[sum, left=0.5cm of a] (sum) {};

\node[left=1cm of a.205] (al){};
\node[right=0.5cm of a.345] (ar){};
\node[left=1.8cm of b] (bl){};
\node[right=1.8cm of b] (br){};
\node[left=1cm of a.155] (kl){};
\node[right= of kl](rkl){};
\node[left=0.5cm of ar] (abr){};
\node[right=1cm of al] (abl){};

\node[bigbox, right=0.55 cm of abr] (ar_PLL) {PLL};
\node[bigbox, left=1.2 cm of abl] (al_CPS) {CPS};
\node[bigbox, above=0.2 cm of al_CPS] (al_d) {Dist};
\node[left=0.5cm of al_d](lkl){$\eta$};
\node[below=0.2cm of lkl](lklt){};
\node[right=0.2cm of lklt](lklb){};

\draw[-] (abr)-- (ar_PLL);
\draw[->] (sum)-- (a)  node[midway,below=0.15 cm]{$i$};

\draw[->] (al_CPS)-|(sum);
\draw[->] (al_d)-|node[near end, right] {$+$} (sum);

\draw[-] (bl.east)-|(lklb);
\draw[->] (lklb)|-(al_CPS.west);

\draw[->] (br)-- (b) node[midway,above]{$y_\text{\it foll}$};
\draw[-] (b)-- (bl.east) node[midway,above]{$u_\text{\it foll}$};

\draw[->] (lkl)-- (al_d);
\draw[-] (2.32,-0.18)-- (2.52,-0.18);
\draw[-] (2.52,-0.173)-- (2.52,-1.44);

\node[right= of a.375](n){};
\node[right=0.7cm of n](p){};
\node[below=0.3cm of p](x){};
\node[right=-0.15cm of x](y){$\underbrace{\begin{bmatrix}\omega_\text{G}\\ \dot{\omega}_\text{G}\\P_\text{G}\\ P_\text{VI} \end{bmatrix}}_{\displaystyle{y_p}}$};
\draw[->] (a.375)-- (p);

\end{tikzpicture}
\caption{Closed-loop system for the grid-following VI with tuning parameters $K_\text{\it foll}$, where CPS is the controllable power source.}
\label{Fig:ClosedSys-Foll}
\end{figure}

For each of the grid-forming devices, we define $y_\text{\it form,k}=(\omega_{\text{VI}, k}, P_{\text{VI}, k})$ as the output, collecting its internal frequency variable and its active power injection; and $u_\text{\it form,k}=\dot{\omega}_{\text{VI}, k}$ as the control input. Choosing $K_{\text{\it form},k}=-\begin{bmatrix}\tilde d_k\, \tilde m_k^{-1}\,\,\tilde m_k^{-1}\end{bmatrix}$, \eqref{eq:VIForm:freq} can be re-expressed as 
\begin{align}
 \dot{\omega}_{\text{VI},k} = u_\text{\it form,k} = K_{\text{\it form},k} \; y_\text{\it form,k}.
\end{align}
The resulting overall system with input $u_\text{\it form}$, output $y_\text{\it form}$, and gain matrix $ K_{\text{\it form}}$ is shown in Figure \ref{Fig:ClosedSys-Form}.
Similarly, for each grid-following device we define $y_\text{\it foll,k}=(\hat{\omega}_k,\dot{\hat{\omega}}_k)$ as the output, collecting the frequency and RoCoF estimates from the PLL; and the active power set-point as the control input $u_\text{\it foll,k}=P^\star_{\text{VI},k}$. With $K_{\text{\it foll},k} =\begin{bmatrix}\tilde{d}_k\,\, \tilde{m}_k \end{bmatrix}$, \eqref{eq:VIFoll} can be re-expressed as 
\begin{align}
 P^\star_{\text{VI},k} = u_\text{\it foll,k} = K_{\text{\it foll},k} \; y_\text{\it foll,k}.
\end{align}
The resulting overall system with input $u_\text{\it foll}$, output $y_\text{\it foll}$, and gain matrix $ K_{\text{\it foll}}$ is shown in Figure \ref{Fig:ClosedSys-Foll}.
Next, let $x$ and $u$ denote the states of the power system equipped with VI devices and the control inputs of these VI devices. Additionally, let $y$ and $y_p$ denote the outputs and the performance outputs corresponding to the performance metrics discussed in Section \ref{subsec:PerfMets}. The overall dynamical model translates to
\begin{subequations}
\label{eq:SysModelCL}
\begin{align}
\dot{x}&=f(x, z, u, \eta), \label{eq:SysModelCL:diff}\\
\vzeros[] &=g(x,z),\label{eq:SysModelCL:alg}\\
(y,\, y_p)&= (h(x),\,h_p(x)).
\end{align}
\end{subequations}
Next, we linearize these dynamics around a nominal operating point. In this process, the algebraic equation \eqref{eq:SysModelCL:alg} can be eliminated by exploiting the fact that its Jacobian with respect to $x$ has full rank at operating points that do not correspond to voltage collapse. Likewise, we can also remove the unobservable mode (with zero eigenvalue) corresponding to absolute angles to obtain a linearization 
\begin{subequations}
\label{eq:Closed-LoopLin}
\begin{align}
\Delta\dot{x} &=A \Delta x + B \Delta u + G \eta,\\
\Delta y &= C \Delta x, \, \Delta y_{p}=C_{p}\Delta x,
\end{align}
\end{subequations}
where  $\Delta x$, $\Delta y$, $\Delta y_{p}$, $\Delta u$ are the resulting deviation states, measurement outputs, performance outputs, control inputs; and ${G}=B\Pi$ is the disturbance gain matrix which encodes (via $\Pi=\text{diag}\{\pi_1, \pi_2, \ldots\}$) the location and (relative) strengths of the disturbances $\eta$. The states $x$ and outputs $y$ are different for both VI implementations. For grid-following implementation, these correspond to $x={\begin{pmatrix} x_s,\, x_\text{PLL} \end{pmatrix}}$, $y_\text{\it foll}={\begin{pmatrix}\hat\omega,\, \dot{\hat\omega} \end{pmatrix}}$ whereas, $x={\begin{pmatrix} x_s,\, x_\text{VI} \end{pmatrix}}$, $y_\text{\it form}={\begin{pmatrix} \omega_\text{VI},\, {P_\text{VI}} \end{pmatrix}}$ refer to grid-forming implementation.

\subsection{Virtual inertia as output feedback}
The control inputs are given by static output feedback, i.e., 
\begin{subequations}
\begin{align}
  \label{eq:ControlFoll}
 \Delta u_\text{\it foll} &= \underbrace{  
  \begin{bmatrix} \tilde{m}_1 &  \tilde{d}_1 & \ldots & 0 & 0 \\
                   \vdots & \vdots  & \ddots & \vdots & \vdots\\  
                     0        & 0 & \ldots & \tilde{m}_{n_c} & \tilde{d}_{n_c} 
 \end{bmatrix}}_{\displaystyle {K_\text{\it foll}}} \underbrace{\begin{bmatrix} \Delta \hat\omega_{1} \\ {\Delta \dot{\hat\omega}_{1}}\\\vdots \\ {\Delta \hat\omega_{n_c}}\\{\Delta \dot{\hat\omega}_{n_c}} \end{bmatrix} }_{\displaystyle {\Delta y_\text{\it foll}}},\\
 \label{eq:ControlForm}
\Delta  u_\text{\it form} &= \underbrace{  
  \begin{bmatrix}  \tilde{\alpha}_1 &  \tilde{\beta}_1 & \ldots & 0 & 0 \\
                   \vdots & \vdots  & \ddots & \vdots & \vdots \\  
                   0          & 0 & \ldots & \tilde{\alpha}_{n_c} & \tilde{\beta}_{n_c} 
  \end{bmatrix}}_{\displaystyle {K_\text{\it form}}} \underbrace{\begin{bmatrix} {\Delta \omega_{\text{VI}, 1}} \\ {\Delta P_{\text{VI}, 1}}\\\vdots \\  {\Delta \omega_{\text{VI}, n_c}}\\{\Delta P_{\text{VI}, n_c}}\end{bmatrix}}_{\displaystyle {\Delta y_\text{\it form}}}, 
\end{align}
\end{subequations}
are the control inputs for the grid-following and grid-forming implementations respectively, with feedback matrices $K_\text{\it foll}$, $K_\text{\it form}$, feedback gains $\tilde{\alpha}_k=-\tilde{d}_k\,{\tilde{m}_k}^{-1} $, $\tilde{\beta}_k={\tilde{m}_k}^{-1}$, and number of virtual inertia devices $n_c$. 

For our analysis, the performance output is selected as
\begin{align*}
  \Delta y_p &= C_p \, \Delta x= \begin{pmatrix}  r^{\frac{1}{2}}_{\omega} \Delta \omega_\text{G},\, r^{\frac{1}{2}}_{\dot \omega}\, \Delta \dot\omega_\text{G},\, r^{\frac{1}{2}}_\text{G} \Delta P_\text{G},\, r^{\frac{1}{2}}_\text{VI} \Delta P_\text{VI} \end{pmatrix},
        \label{eq:perfOut}
\end{align*}
where the states $x$ depend on the VI implementation. For impulse disturbances, the infinite horizon integral of quadratic penalties on frequency deviations, RoCoF, as well as power injections from VI devices and generators is given by
\begin{align}
J_{\infty, \delta(t)}= \int_{0}^{\infty} \Delta  y_p^{\top} \Delta  y_p \, \text{d}t.
\end{align}
With $A_\text{\it cl}=A+BKC$, and combining \eqref{eq:Closed-LoopLin} and \eqref{eq:ControlFoll} (respectively, \eqref{eq:ControlForm}) the resulting dynamical system $\mathcal G$ is 
\begin{align}
 \Delta {\dot{x}} = A_\text{\it cl} \Delta x + {G} {\eta}, \quad 
 \Delta  y_p = {C_p}\, \Delta {x}\, .\label{eq:sysH2cl}
\end{align}

\subsection{$\mathcal H_{2}$ norm optimization}
To compute the $\mathcal H_{2}$-norm between the disturbance input ${\eta}$ and the performance output $y_p$ of the system \eqref{eq:sysH2cl}, let the so-called observability Gramian ${P}_K$ denote the positive definite solution $P$ of the Lyapunov equation
\begin{align}
\label{eq:h2lyap}
 {P} A_\text{\it cl}+A_\text{\it cl}^{\top} {P}  + C_p^{\top} {C_p} = \vzeros[],
\end{align}
parameterized in ${K}$ for the given system matrices ${A}$, ${B}$, ${C}$, and ${C_p}$. The $\mathcal H_2$ norm is given by $J_{\infty, \delta (t)} = \|\mathcal G\|^2_2 = \trace({G}^\top {P}_{K} {G})$ (cf. \cite{KZ-JCD-KG:96}). Thus, the optimization problem to compute the $\mathcal H_2$ optimal allocation is 
\begin{align}
\label{eq:h2opt}
 \min_{K \in \mathcal{S} \cap \, \mathcal{C}} \quad J_{\infty, \delta (t)}.
\end{align}
The set $\mathcal{S}$ is used to encode the structural constraint on ${K}$, i.e., the purely local feedback structure of the virtual inertia control in \eqref{eq:ControlFoll} and \eqref{eq:ControlForm}. Hereafter, $\mathcal{C}$ denotes the set of constraints on the control gains discussed in Section \ref{sec:constraint}. 

The optimization problem \eqref{eq:h2opt} can tune the gain of any VI device in the system. Sparse allocations (i.e., with few VI devices with significant contribution) can be obtained by including an $\ell_1$-penalty in the optimization \cite[Sec. 3.5]{BKP-SB-FD:17}. Note that evaluating the cost function requires solving the Lyapunov equation \eqref{eq:h2lyap}, which is non-linear in $P$ and $K$. In general, the optimization problem \eqref{eq:h2opt} is non-convex and may be of very large-scale. However, by exploiting the feedback structure of the problem, the gradient of  $\|\mathcal G\|^2_2$ with respect to $K$ can be computed efficiently and can be directly used to solve \eqref{eq:h2opt} via scalable first order methods (e.g., projected gradient) or to speed up higher order methods (see the Appendix for details). 

\subsection{Complexity of the gradient computation}
In \cite{BKP-SB-FD:17} gradient-based optimization methods are used to directly optimize the inertia constants of a linearized networked swing equation model to minimize the $\mathcal H_2$ norm of a power system. For a system with $n$ buses, the gradient computation in \cite{BKP-SB-FD:17} requires the solution of $n-1$ Lyapunov equations of dimension $2n$, resulting in a complexity of $\mathcal{O}((n+1)n^3)$. In contrast \eqref{eq:h2opt} includes more realistic models of virtual inertia devices and the gradient of $J_{\infty, \delta (t)}$ can be computed by solving two Lyapunov equations of dimension $4n$, thereby reducing the complexity to $\mathcal{O}(n^3)$.  In \cite{TSB-FD:17} a sequential linear programming approach is used to optimize the allocation of grid-following virtual inertia and damping. This method directly optimizes the frequency nadir and RoCoF. However, every iteration of the optimization algorithm in \cite{TSB-FD:17} requires computing the eigenvalues of the linearized system which has complexity $\mathcal{O}(n^3)$, as well as time-domain simulations and the solution of a linear program, resulting in far higher computation complexity than the proposed method.

\section{Test case Description}\label{sec:test}
\begin{figure}[hbp]
\centering
\includegraphics[height=1.5\columnwidth]{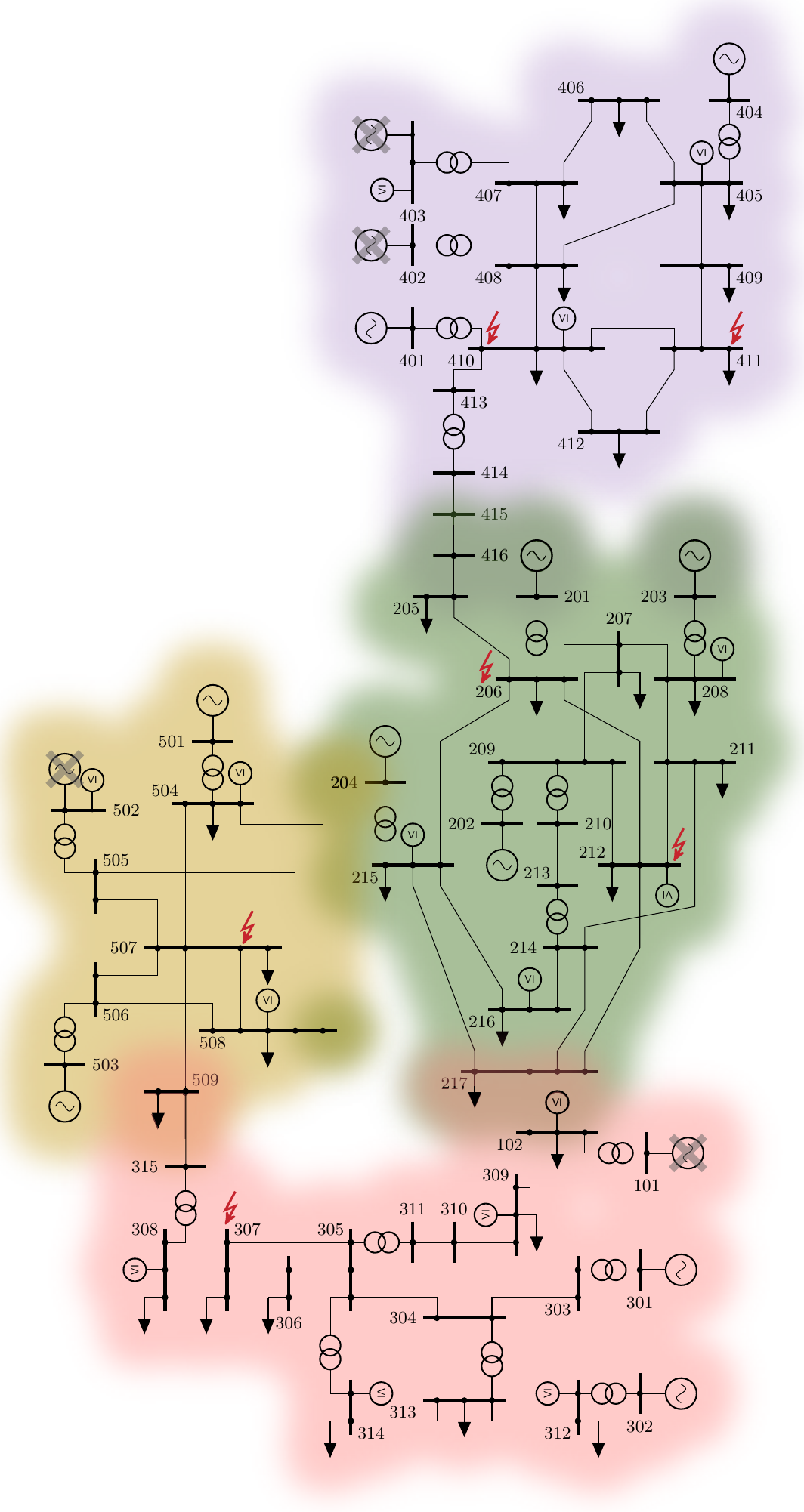}
\caption{South-East Australian Power System line diagram. The crossed out generators are replaced by constant power sources to mimic a low-inertia scenario, whereas the red lightning symbols are the locations where disturbances are injected. The circles with VI inscribed within indicate the virtual inertia and damping devices distributed across the power system.}
\label{Fig:AusGrid}
\end{figure}

To illustrate our algorithms for optimal inertia and damping tuning, we use a test case based on the 14-generator, 59-bus South-East Australian system \cite{MG-DV:14,ASA-SR-GV-AC-DJH:17} shown in Figure~\ref{Fig:AusGrid}. 

It is equipped with higher order models for turbines, governors, power system stabilizers (PSSs), and voltage regulators (AVRs). This system has several interesting features, for instance its string topology and weak coupling between South Australia (area 5) and the rest of the system. The SIMULINK version \cite{AM-IK-PB-GS:15} of the model \cite{MG-DV:14} was developed for the light loading scenario. Variations of this model have also been studied as low-inertia test cases in \cite{TSB-FD:17, SPL-PM:18}.

For this paper the model presented in \cite{AM-IK-PB-GS:15} was modified to obtain a low-inertia case study by replacing synchronous machines located at the buses labeled 101, 402, 403, and 502 with constant power sources\footnote{In other words, sources that maintain constant active and reactive power injections regardless of the frequency or voltage at their point of connection.} that inject the same active and reactive power as the original generators. This modeling choice is based on the high penetration of renewable generation in the real-world power system (particularly in area 5) that does not provide frequency support \cite{SPL-PM:18}. The model was augmented with 15 VI devices across the system (see Figure~\ref{Fig:AusGrid}). For brevity of the presentation we consider two scenarios, in the first scenario the VI devices are all grid-forming, in the second scenario they are all grid-following (see Section~\ref{sec:model}). In the case study in \cite{TSB-FD:17} motor loads with non-negligible inertia are used to ensure that the notion of a frequency signal (as input the VI devices) is well defined. In this work, we do not require this assumption. Finally, we use constant power injections at six locations (indicated by a red bolt) to simulate disturbances. The SIMULINK model of the benchmark system including virtual inertia devices is available online \cite{BKP-DG:18}.

\section{Results}
\label{sec:results}
In this section we compare the performance of the original system with the closed-loop system equipped with virtual inertia and damping devices. We consider both the grid-following and the grid-forming modes of implementation and mainly focus on the performance metrics defined in Section~\ref{subsec:PerfMets}.

\subsection{Validity of the linearized model}
As discussed in Section \ref{sec:metrics}, we optimize the virtual inertia and damping gains using a linearization of the system at the nominal operating point. To validate the linearized model we compare it to the non-linear model for step disturbances at the six locations shown in Figure~\ref{Fig:AusGrid} ranging from $-250\, \mathrm{MW}$ to $+250\, \mathrm{MW}$. In Figure~\ref{fig:lin}, the relative linearization errors for different performance metrics are plotted- both for the grid-following and the grid-forming virtual inertia and damping implementations. The plots reveal a concentration of data points in the $-10\%$ to $+10\%$ band. This indicates that the linear approximation of the model closely resembles the non-linear model and justifies the effectiveness of our approach. 

\subsection{Optimal tuning and placement of VI devices}
The optimal inertia and damping profiles for the system are computed using the optimization problem \eqref{eq:h2opt}. We consider the same weighted performance outputs \eqref{eq:perfOut} for both grid-forming and grid-following and set the penalties to $r_{\omega}=0.1$, $r_{\dot{\omega}}\!=\!0.2$, $r_\text{G}\!=\!0.2$, and $r_\text{VI}\!=\!0.2$, thereby identically penalizing the power injections from the VI devices and the synchronous machines. Further, the disturbance gain matrix $\Pi$, introduced in Section \ref{sec:clandlin} is set to identity, i.e., $\pi_i=1$. In other words each node is subject to equally sized disturbances. Finally, using the approach outlined in Section \ref{sec:constraint} we choose the constraints such that $\sum_k \tilde{d}_k \leq 420\, \mathrm{MW s/rad}$,  $\tilde{d}_k \leq 40 \,\mathrm{MW s/rad}$, and $\tilde{m}_k \leq 18.5\, \mathrm{MW s^2/rad}$. 
These constraints ensure that the total damping does not exceed realistic values, and that the power output of the converters is roughly limited to $40\, \mathrm{MW}$ for frequency deviations in the normal operating regime. 
The resulting inertia and damping allocations for the above parameters and constraints are depicted in Figure~\ref{fig:allocs} {\it (a)}, {\it (b)}. 
\begin{figure}[t]
%
%
\definecolor{mycolor1}{HTML}{F1931B}%
\begin{tikzpicture}

\begin{axis}[%
width=1.25in,
height=0.75in,
at={(0.683in,1.945in)},
scale only axis,
xmin=-50,
xmax=50,
xtick={ -50, -25, 0, 25, 50},
xlabel style={font=\color{black}},
xlabel={\small $\Delta \left\vert{\omega}_{\text{G},k} \right\vert_\infty$ [\%]},
ymin=0,
ymax=0.35,
ytick={ 0, 0.1, 0.2, 0.3},
yticklabels={$0$,$10$,$20$,$30$},ylabel style={font=\color{black}},
ylabel={\small samples [\%]},
axis background/.style={fill=white},
xmajorgrids,
ymajorgrids,
yticklabel style = {font=\footnotesize,xshift=0ex},
xticklabel style = {font=\footnotesize,yshift=0ex},
every axis plot/.append style={ultra thin},
legend style={legend cell align=left, align=left, draw=white!15!black}
]

\addplot[ybar interval, fill=mycolor1, fill opacity=0.9, area legend] table[row sep=crcr] {%
x	y\\
-85	0.00222222222222222\\
-80	0\\
-75	0.00333333333333333\\
-70	0.00777777777777778\\
-65	0.00555555555555556\\
-60	0\\
-55	0.00666666666666667\\
-50	0\\
-45	0.00111111111111111\\
-40	0.00111111111111111\\
-35	0.00555555555555556\\
-30	0.0166666666666667\\
-25	0.00555555555555556\\
-20	0.00333333333333333\\
-15	0.0433333333333333\\
-10	0.101111111111111\\
-5	0.284444444444444\\
0	0.307777777777778\\
5	0.133333333333333\\
10	0.00555555555555556\\
15	0.00888888888888889\\
20	0.0277777777777778\\
25	0.00666666666666667\\
30	0\\
35	0.00222222222222222\\
40	0.0122222222222222\\
45	0.00777777777777778\\
50	0.00777777777777778\\
};
\end{axis}

\begin{axis}[%
width=1.25in,
height=0.75in,
at={(2.458in,1.945in)},
scale only axis,
xmin=-50,
xmax=50,
xtick={ -50, -25, 0, 25, 50},
xlabel style={font=\color{black}},
xlabel={\small $\Delta \left\vert\dot\omega_{\text{G},k}\right\vert_{\infty}$ [\%]},
ymin=0,
ymax=0.35,
ytick={  0, 0.1, 0.2, 0.3},
yticklabels={$0$,$10$,$20$,$30$},ylabel style={font=\color{black}},
ylabel={\small samples [\%]},
axis background/.style={fill=white},
xmajorgrids,
ymajorgrids,
yticklabel style = {font=\footnotesize,xshift=0ex},
xticklabel style = {font=\footnotesize,yshift=0ex},
every axis plot/.append style={ultra thin},
legend style={legend cell align=left, align=left, draw=white!15!black}
]

\addplot[ybar interval, fill=mycolor1, fill opacity=0.9, area legend] table[row sep=crcr] {%
x	y\\
-60	0.00138888888888889\\
-55	0.00277777777777778\\
-50	0.00138888888888889\\
-45	0.00277777777777778\\
-40	0.00694444444444444\\
-35	0.00277777777777778\\
-30	0.00416666666666667\\
-25	0.00972222222222222\\
-20	0.00694444444444444\\
-15	0.0333333333333333\\
-10	0.105555555555556\\
-5	0.3125\\
0	0.331944444444444\\
5	0.109722222222222\\
10	0.0236111111111111\\
15	0.0125\\
20	0.00972222222222222\\
25	0.0111111111111111\\
30	0\\
35	0\\
40	0.00277777777777778\\
45	0.00694444444444444\\
50	0\\
55	0\\
60	0.00138888888888889\\
65	0.00138888888888889\\
};

\end{axis}

\begin{axis}[%
width=1.25in,
height=0.75in,
at={(0.683in,0.709in)},
scale only axis,
xmin=-50,
xmax=50,
xtick={ -50, -25, 0, 25, 50},
xlabel style={font=\color{black}},
xlabel={\small $\Delta \left\vert P_{\text{VI},k} \right\vert_{\infty}$ [\%]},
ymin=0,
ymax=0.35,
ytick={  0, 0.1, 0.2, 0.3},
yticklabels={$0$,$10$,$20$,$30$},ylabel style={font=\color{black}},
ylabel={\small samples [\%]},
axis background/.style={fill=white},
xmajorgrids,
ymajorgrids,
yticklabel style = {font=\footnotesize,xshift=0ex},
xticklabel style = {font=\footnotesize,yshift=0ex},
every axis plot/.append style={ultra thin},
legend style={legend cell align=left, align=left, draw=white!15!black}
]

\addplot[ybar interval, fill=mycolor1, fill opacity=0.9, area legend] table[row sep=crcr] {%
x	y\\
-85	0.00138888888888889\\
-80	0\\
-75	0\\
-70	0.00694444444444444\\
-65	0.00277777777777778\\
-60	0.00694444444444444\\
-55	0.00277777777777778\\
-50	0.00138888888888889\\
-45	0.00138888888888889\\
-40	0\\
-35	0.00694444444444444\\
-30	0.0152777777777778\\
-25	0.00277777777777778\\
-20	0.00138888888888889\\
-15	0.0430555555555556\\
-10	0.0986111111111111\\
-5	0.305555555555556\\
0	0.308333333333333\\
5	0.141666666666667\\
10	0.00277777777777778\\
15	0.0111111111111111\\
20	0.0138888888888889\\
25	0\\
30	0.0111111111111111\\
35	0\\
40	0.00833333333333333\\
45	0.00555555555555556\\
50	0.00555555555555556\\
};


\end{axis}

\begin{axis}[%
width=1.25in,
height=0.75in,
at={(2.458in,0.709in)},
scale only axis,
xmin=-50,
xmax=50,
xtick={ -50, -25, 0, 25, 50},
xlabel style={font=\color{black}},
xlabel={\small $\Delta \left\vert P_{\text{G},k} \right\vert_{\infty}$ [\%]},
ymin=0,
ymax=0.35,
ytick={  0, 0.1, 0.2, 0.3},
yticklabels={$0$,$10$,$20$,$30$},ylabel style={font=\color{black}},
ylabel={\small samples [\%]},
axis background/.style={fill=white},
xmajorgrids,
ymajorgrids,
yticklabel style = {font=\footnotesize,xshift=0ex},
xticklabel style = {font=\footnotesize,yshift=0ex},
every axis plot/.append style={ultra thin},
legend style={legend cell align=left, align=left, draw=white!15!black}
]

\addplot[ybar interval, fill=mycolor1, fill opacity=0.9, area legend] table[row sep=crcr] {%
x	y\\
-85	0.00111111111111111\\
-80	0\\
-75	0\\
-70	0.00666666666666667\\
-65	0.00333333333333333\\
-60	0\\
-55	0.00222222222222222\\
-50	0.00666666666666667\\
-45	0.00111111111111111\\
-40	0\\
-35	0.00555555555555556\\
-30	0.01\\
-25	0.00777777777777778\\
-20	0.00222222222222222\\
-15	0.0433333333333333\\
-10	0.0933333333333333\\
-5	0.314444444444444\\
0	0.304444444444444\\
5	0.147777777777778\\
10	0.00111111111111111\\
15	0.00888888888888889\\
20	0.0177777777777778\\
25	0.00777777777777778\\
30	0.00111111111111111\\
35	0\\
40	0.00888888888888889\\
45	0.00444444444444444\\
50	0.00444444444444444\\
};


\end{axis}
\end{tikzpicture}%
\caption{Distribution of relative linearization errors for load steps ranging from $-250\, \mathrm{MW}$ to $+250\, \mathrm{MW}$ at the nodes indicated in Figure \ref{Fig:AusGrid} for both grid-forming, grid-following configurations.\label{fig:lin}}
\end{figure}
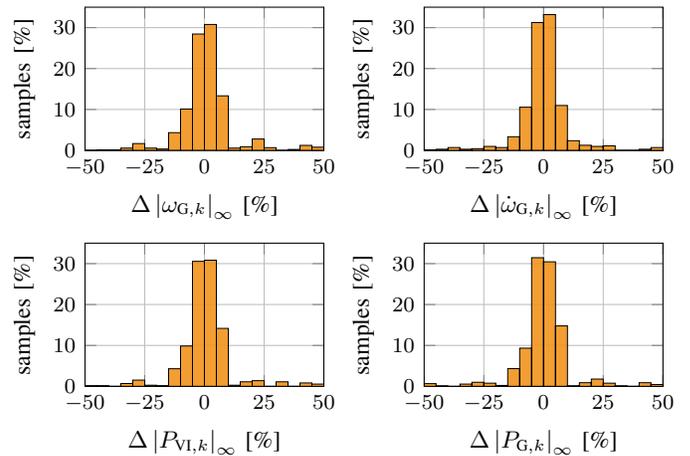
\begin{figure}[b!!]
%
%
\definecolor{mycolor1}{HTML}{A4978E}%
\definecolor{mycolor3}{HTML}{BE9063}%
\definecolor{mycolor2}{HTML}{525B56}%
\definecolor{mycolor4}{HTML}{132226}%

\vspace{-1em}
\begin{tikzpicture}

\begin{axis}[%
width=3.1in,
height=1.2in,
at={(1.011in,2.80in)},
scale only axis,
bar shift auto,
xmin=0.514285714285714,
xmax=15.4857142857143,
xtick={1,2,3,4,5,6,7,8,9,10,11,12,13,14,15},
xticklabels={{102},{208},{212},{215},{216},{308},{309},{312},{314},{403},{405},{410},{502},{504},{508}},
ymin=0,
ymax=58,
ytick={0, 10, 20, 30, 40, 50},
yticklabel style = {font=\footnotesize,xshift=0ex},
xticklabel style = {font=\footnotesize,yshift=0ex},
ymajorgrids,
xlabel={\small node},
every axis plot/.append style={ultra thin},
axis background/.style={fill=white},
title style={font=\it},
title={\small{(a) Grid-Forming}},
legend style={at={(0.03,0.97)}, anchor=north west, legend columns=-1, legend cell align=left, align=left, draw=none}
]
\addplot[ybar, bar width=4, fill=mycolor3, draw=black, area legend] table[row sep=crcr] {%
1	4.75198149650708\\
2	4.80203238758951\\
3	4.84224607047909\\
4	4.79996255974154\\
5	4.84951678090365\\
6	4.81797857617476\\
7	4.67207909454179\\
8	2.9201776381046\\
9	4.79870163095871\\
10	3.90289585465935\\
11	3.65582194279329\\
12	2.80930569045934\\
13	17.9748721616244\\
14	16.5146163338333\\
15	13.7552983390899\\
};
\addplot[forget plot, color=white!15!black] table[row sep=crcr] {%
0.514285714285714	0\\
15.4857142857143	0\\
};
\addlegendentry{\small{$\text{inertia [MW s}^\text{2}\text{/rad]}\quad$}}

\addplot[ybar, bar width=4, fill=mycolor4, draw=black, area legend] table[row sep=crcr] {%
1	23.2801623293638\\
2	35.8974828213254\\
3	37.0652403971668\\
4	34.565880334748\\
5	32.7234704046107\\
6	37.5254758727884\\
7	24.4469532756254\\
8	5.07155242999244\\
9	37.984191746052\\
10	14.580946382229\\
11	9.64214118054151\\
12	6.38799973210701\\
13	27.1065997687833\\
14	29.5731985809343\\
15	20.0522664371918\\};
\addplot[forget plot, color=white!15!black] table[row sep=crcr] {%
0.514285714285714	0\\
15.4857142857143	0\\
};
\addlegendentry{\small{damping [MW s/rad]}}
\end{axis}

\begin{axis}[%
width=3.1in,
height=1.2in,
at={(1.011in,0.9in)},
scale only axis,
bar shift auto,
xmin=0.514285714285714,
xmax=15.4857142857143,
xtick={1,2,3,4,5,6,7,8,9,10,11,12,13,14,15},
xticklabels={{102},{208},{212},{215},{216},{308},{309},{312},{314},{403},{405},{410},{502},{504},{508}},
ymin=0,
ymax=58,
ytick={0, 10, 20, 30, 40, 50},
yticklabel style = {font=\footnotesize,xshift=0ex},
xticklabel style = {font=\footnotesize,yshift=0ex},
xlabel={\small node},
ymajorgrids,
axis background/.style={fill=white},
every axis plot/.append style={ultra thin},
title style={font=\it},
title={\small{(b) Grid-Following}},
legend style={at={(0.03,0.97)}, anchor=north west, legend columns=-1, legend cell align=left, align=left, draw=none}
]
\addplot[ybar, bar width=4, fill=mycolor3, draw=black, area legend] table[row sep=crcr] {%
1	1.53626583555718\\
2	3.41260616749131\\
3	0.487700184178425\\
4	7.26850155138847\\
5	2.50402482335262\\
6	0.198429745379375\\
7	0.928559523627684\\
8	14.2079110661918\\
9	0.595669598174814\\
10	1.22022251218213\\
11	5.45076812688971\\
12	18.5\\
13	18.5\\
14	18.5\\
15	18.5\\
};
\addplot[forget plot, color=white!15!black] table[row sep=crcr] {%
0.514285714285714	0\\
15.4857142857143	0\\
};
\addlegendentry{\small{$\text{inertia [MW s}^\text{2}\text{/rad]} \quad$}}

\addplot[ybar, bar width=4, fill=mycolor4, draw=black, area legend] table[row sep=crcr] {%
11	7.79015816226344\\
2	40\\
3	40\\
4	37.3386452070442\\
5	32.3090753427562\\
6	12.8645833166182\\
7	5.19055009628645\\
8	40\\
9	29.4575569562099\\
10	-0\\
11	15.0494309188216\\
12	40\\
13	40\\
14	40\\
15	40\\
};
\addplot[forget plot, color=white!15!black] table[row sep=crcr] {%
0.514285714285714	0\\
15.4857142857143	0\\
};
\addlegendentry{\small{damping [MW s/rad]}}

\end{axis}
\end{tikzpicture}%
\caption{Optimal inertia and damping allocations for the Australian system for the grid-forming and grid-following configurations.\label{fig:allocs}}
\end{figure}
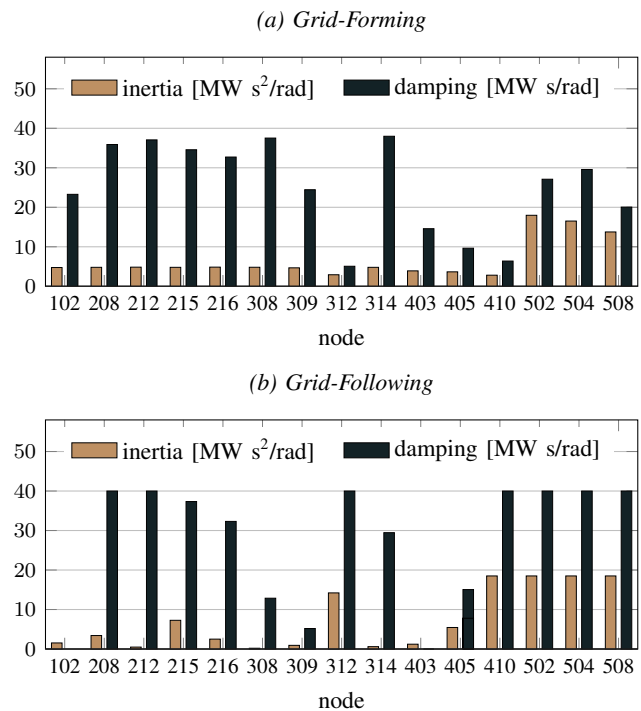
Finally, we observe that no significant performance gains can be achieved by optimizing the PLL gains beyond applying standard tuning techniques (see \cite{C00,GMF+14}).

\subsection{Contrasting allocations for different VI implementations}
\begin{figure}[htp]
\input{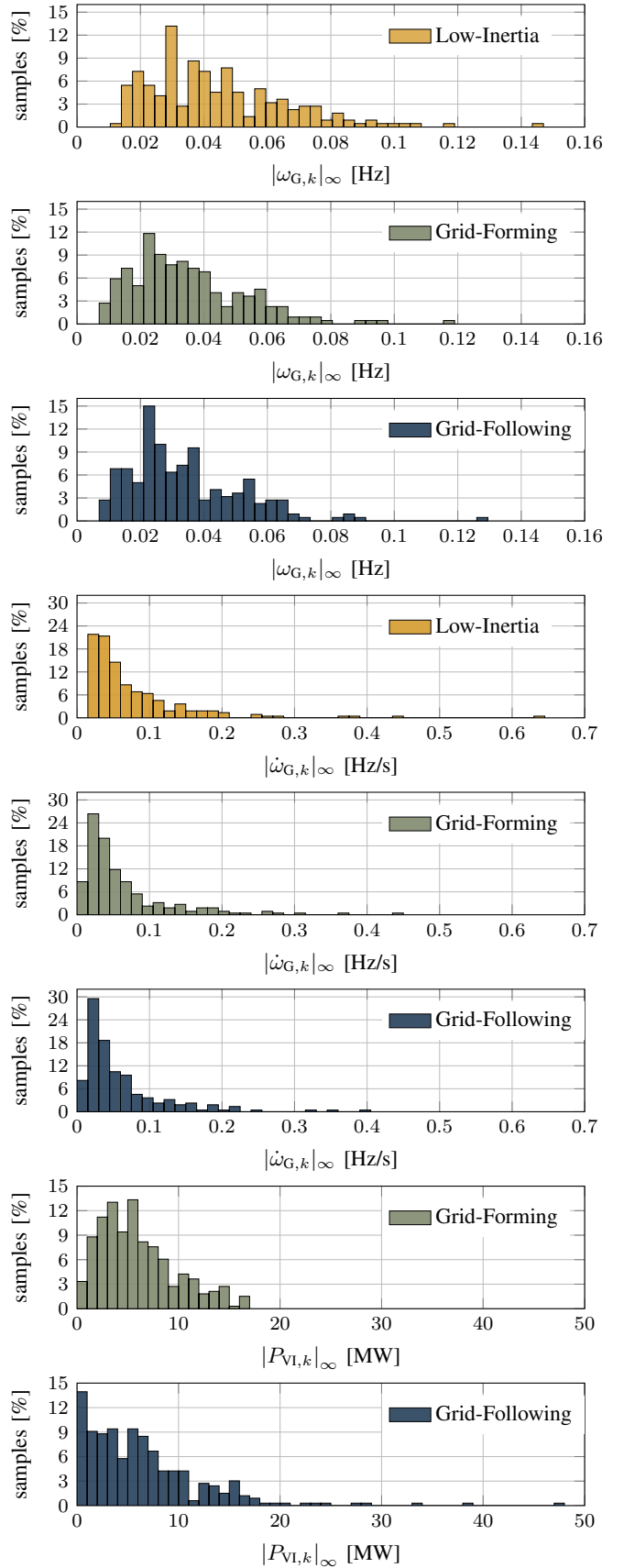}
 \caption{Distribution of generator frequency nadirs, maximum generator RoCoF, and VI power injections for load steps ranging from $-350 \,\mathrm{MW}$ to $-150\, \mathrm{MW}$ at the nodes indicated in Figure \ref{Fig:AusGrid} for different converter configurations.}
 \label{fig:hist_randomised}
\end{figure}
The two allocations highlight some interesting features. Note that the optimized allocations are not uniform across the system.  In fact, the virtual inertia for both implementations is predominantly allocated in area 5. Incidentally, the blackout in South-Australia in 2016 was also in this area \cite{AEMO:17}. Moreover, uniform allocations, chosen as the initial guess for the optimization, are typically not optimal (see also \cite{BKP-SB-FD:17}). Another facet of the allocations is that the gains for the grid-following virtual inertia devices are limited by the constraints imposed in the optimization. This may be primarily attributed to time-delays (RoCoF estimation, response time $\tau_\text{\it foll}$ of the power source, etc.) encountered for the inertial response. To compensate for these delays, the allocation for the grid-following VI relies on larger inertia and damping gains (see Figures~\ref{fig:allocs}) at some nodes as well as larger total damping and inertia (refer Table~\ref{tab:case}). While significant inertia and damping is allocated at all nodes in the case of grid-forming virtual inertia devices, the grid-following implementation results in negligible allocations for some nodes outside of area 5.

\subsection{Impact of VI devices on frequency stability}
To investigate the effect of the VI devices, we consider the non-linear model of the South East Australian grid with the optimal VI device parameter values from Figure~\ref{fig:allocs}. Next, individual step disturbances at the six locations, as shown in Figure~\ref{Fig:AusGrid} are considered. These disturbances range from $-375 \,\mathrm{MW}$ to $-150 \,\mathrm{MW}$ and capture a load increase (or equivalently a loss of generation). In Figure~\ref{fig:hist_randomised}, the resulting distribution of post-fault frequency nadirs {$\left\vert{\omega}_{\text{G},k} \right\vert_\infty$}, maximum RoCoF {$\left\vert{\dot\omega}_{\text{G},k} \right\vert_\infty$} of all generators for the system without VI devices, the system with grid-following VI devices, and the system with grid-forming VI devices are represented along with the peak power injections {$\left\vert P_{\text{VI},k} \right\vert_{\infty}$} from the VI devices. We make the following observations:
\begin{enumerate}[label=(\alph*), leftmargin=0.1cm, itemindent=0.5cm]
\setlength{\itemsep}{2pt}
\item The mean and variance of the distribution of frequency nadirs is smaller for the system equipped with the VI devices as compared to system without VI devices. 
\item The grid-following implementation has a smaller mean frequency nadir, but larger variance and longer tail in comparison to grid-forming counterpart. However, this comes at the expense of a much larger peak power injection for certain disturbances. Note that, the maximum injection for the grid-following VI, for the same set of disturbances is roughly three times that of the grid-forming VI implementation. 
\item While we observe a smaller mean and a shorter tail of the distribution of the maximum generator RoCoF for grid-following virtual inertia, the histograms suggest that the impact of the virtual inertia devices on the maximum RoCoF is modest. However, it is noteworthy that these histograms only depict the maximum RoCoF {$\left\vert {\dot\omega}_{\text{G},k} \right\vert_\infty$} at generator buses. Moreover, the maximum is typically attained during the first swing of the system after a fault. In contrast, the time-domain simulations depicted in Figure~\ref{fig:time-dom} show that virtual inertia devices can have significant impact on the RoCoF $\dot{\omega}_\text{G}$ after the first swing. These differences are not captured when using the maximum RoCoF {$\left\vert {\dot\omega}_{\text{G},k} \right\vert_\infty$} as performance metric, but are accurately captured by the $\mathcal{H}_2$ norm.
\end{enumerate}
We conclude that the VI devices have the expected positive impact on frequency stability. Moreover, the differences between the two VI implementations appear to be mostly related to differences in the maximum power injection. In the next section we will investigate the time-domain response of the system with and without VI devices in more detail.

\subsection{Time-domain responses}\label{subsec:time-dom_resp}
\begin{figure}[b!!]
\input{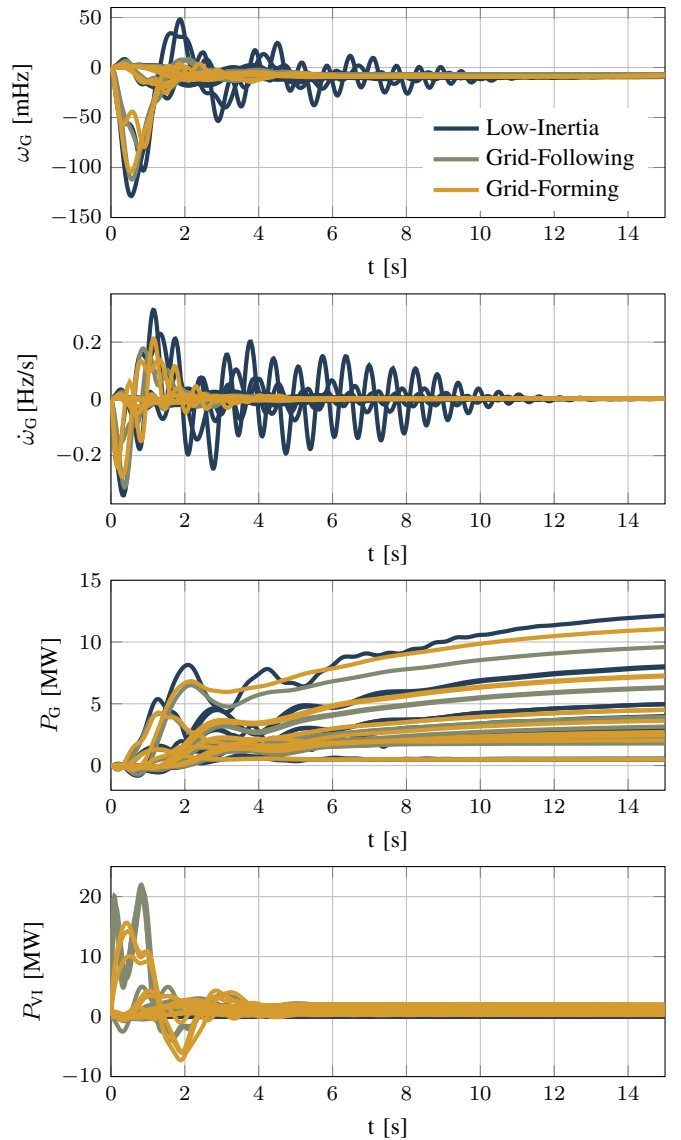}
 \caption{Time-domain plots for generator frequencies, generator RoCoF, generator power injections, and power injections of the VI devices for the low-inertia  system, grid-following, and grid-forming configurations for a step disturbance of $200 \, \mathrm{MW}$ at node 508.}
 \label{fig:time-dom}
\end{figure}

\begin{table}[b!!!]
\begin{center}
\caption{Performance metrics for a $200\, \mathrm{MW}$ load increase at node 508}
\label{tab:case}
\bgroup
\def\arraystretch{1.35}
\setlength{\tabcolsep}{1.5pt}
\begin{tabular}{p{3.30cm} R{1.30cm} R{1.80cm} R{1.80cm}}
\hline\noalign{\smallskip}
\multicolumn{1}{l}{Performance metric} &  {Original}   & {Grid-Following} & {Grid-Forming} \\
\noalign{\smallskip}\hline\noalign{\smallskip}
$\sum_i \tilde{m}_i$ [$\mathrm{{MW}\,{s^2}/rad}$] & - & $111.8$ & $99.9$ \\
$\sum_i \tilde{d}_i$ [$\mathrm{{MW}\,{s}/rad}$] & - & $420$ & $375.9$ \\
$\max_k \left\vert\dot\omega_{\text{G},k}\right\vert_{\infty}$ [$\mathrm{{Hz}/{s}}$] & $0.34$ & $0.31$ & $0.27$  \\
$\max_k \left\vert\omega_{\text{G},k}\right\vert_{\infty}$ [$\mathrm{mHz}$] & $128.6$ & $112.1$ & $104.3$ \\
$\max_k \left\vert P_{\text{VI},k} \right\vert_{\infty}$  [$\mathrm{MW}$]  & - & $21.98$ & $15.62$  \\
$\max\limits_{t \geq 0} \left\vert \sum_k {P}_{\mathrm{VI},k}(t)\right\vert$  [$\mathrm{MW}$]  & - & $68.1$ & $41.1$  \\
$\max\limits_{t \geq 0} \left\vert \sum_k {P}_{\mathrm{G},k}(t)\right\vert$  [$\mathrm{MW}$]  & $50.5$ &  $39.9$ & $45.9$ \\
$\mathcal{H}_2$ norm & $11.72$ & $9.92$ & $9.66$\\
\noalign{\smallskip}\hline\noalign{\smallskip}
\end{tabular}
\egroup
\end{center}
\end{table}
We now simulate a load increase of $200\, \mathrm{MW}$ at node $508$, this represents a realistic contingency in the system (see \cite[Sec. II]{ASA-SR-GV-AC-DJH:17}). Broadly speaking, this disturbance could also represent a loss of $200\, \mathrm{MW}$ renewable generation in area 5 and is of the type considered in the $\mathcal{H}_2$ optimization. Moreover, due to the low levels of rotational inertia in area 5, the placement of this fault corresponds to the worst-case location. The system responses are illustrated in Figure~\ref{fig:time-dom} and underscore the efficacy of virtual inertia and damping devices in a low-inertia power system. The grid-following and grid-forming VI implementation with the optimal allocations from Figure~\ref{fig:allocs} are simulated and compared with the response of the low-inertia system. Table~\ref{tab:case} shows the key performance indicators discussed in Section~\ref{subsec:PerfMets}. The top two panels of the time-domain plots in Figure~\ref{fig:time-dom} illustrate the frequencies and the RoCoF of the 10 generators in the system for the two different VI implementations and the low-inertia system. The power injections from the generators and the 15 virtual inertia and damping devices across the power system are plotted in the bottom two panels of the figure. The key insights drawn from a closer analysis of these plots are summarized below:

\begin{enumerate}[label=(\alph*), leftmargin=0.1cm, itemindent=0.5cm]
\setlength{\itemsep}{2pt}
\item While both VI implementations improve the frequency nadir and maximum  RoCoF, the grid-forming VI implementation performs better in terms of the absolute values. Further, the total inertia and damping is also less.

\item The maximum active power $\max_k \left\vert P_{\text{VI},k}\right\vert$, injected by a single virtual inertia device as well as the maximum power $\max_{t \geq 0} \left\vert \sum_k {P}_{\mathrm{VI},k}(t)\right\vert$ injected by all the virtual inertia devices combined is smaller for the grid-forming virtual inertia devices. Thus, grid-forming virtual inertia achieves a better performance with a lesser control effort in comparison to the grid-following converters.

\item A drop in the maximum governor response $\max_k \left\vert P_{\text{G},k}\right\vert$, by a single synchronous machine compared to the low-inertia system is observed due to the active power injections from the virtual inertia devices. 

\item A decrease in the $\mathcal{H}_2$ norm is observed for both VI implementations, i.e., the $\mathcal{H}_2$ norm is an effective proxy for time-domain metrics for power system analysis \cite{DG-SB-BKP-FD:17}.

\item Another difference in the implementations pertains to the computation times for solving the optimization problem. Using MATLAB on a Core i7-6600U CPU, the optimization for grid-forming VI takes around $60$s in comparison to $160$s for the grid-following VI for identical penalties.
\end{enumerate}

\section{Summary and Conclusions}
\label{sec:conclusions}
In this paper we considered the problem of low-inertia power systems equipped with grid-following or grid-forming VI implementations using power electronic interfaced renewable energy sources. We modeled these two implementations as dynamic feedback control loops that provide virtual inertia and damping. A system norm-based optimization approach was used to study the problem of optimal placement and tuning of these devices. Our proposed tuning algorithm was far more scalable and computationally efficient in comparison to some of the other existing approaches in the literature. Further, we showcased the capabilities of such VI devices on a high-fidelity non-linear model of the South-East Australian power system and illustrated their efficacy. For a range of disturbances both types of virtual inertia implementations improved the system resilience compared to the system without virtual inertia. The results show that the system robustness does not only depend on the amount of virtual inertia used but can also depend on the specific implementation and location of virtual inertia. However, this fact is in contrast to typical paradigms of ancillary service markets that value energy or the total amount of damping and inertia but not location or specifics of the implementation, and therefore do not capture this aspect. A preliminary result on a market mechanism that considers the location of virtual inertia devices can be found in \cite{BP-SB-LN-FD:17}. Given that our proposed tuning algorithm is computationally efficient, our approach can be used to optimize a virtual inertia allocation with respect to multiple linearized models, each modeling different dispatch points and changes in model structure (e.g., system splits, tripping of generators, etc.). Finally, in future systems operating entirely based on converter-interfaced generation further services that are provided by synchronous machines today (e.g., voltage regulation) need to be provided by grid-forming power converters. Therefore, an interesting direction for future research would involve extending the proposed framework by incorporating suitable performance metrics for ancillary services apart from inertia and fast frequency response.

\section{Acknowledgements}
The authors wish to thank Saverio Bolognani, Theodor Borsche, and Damian Flynn for their fruitful comments.

\appendix
\label{app:grad}
 By using the implicit linearization technique from \cite{TR-EWS:97}, the gradient of the norm $\|\mathcal G\|^2_2$ with respect to ${K}$ is given by
\begin{align}\label{eq:H2grad}
 \nabla_{{K}} \|\mathcal G\|^2_2 &=  2 ({B}^\top {P}_K) {L}_K C^\top,
\end{align}
where ${L}_K$ is the positive semidefinite controllability Gramian obtained as a solution $L$ to the Lyapunov equation
\begin{align}\label{eq:h2lyapL}
{L} {A_\text{\it cl}}^\top + A_\text{\it cl}{L}  + {G}{G}^\top = \vzeros[]\,,
\end{align}
and parameterized in ${K}$ for the given system matrices ${A}$, ${B}$, ${C}$, and ${G}$. Thus, computing the norm
$\|\mathcal G\|^2_2$ and its gradient $\nabla_{K} \|\mathcal G\|^2_2$ for a given $K$ requires solving the Lyapunov equations \eqref{eq:h2lyap} and \eqref{eq:h2lyapL}. Moreover, the number of decision variables of the optimization problem \eqref{eq:h2opt} can be reduced by projecting the gradient 
$\nabla_{K} \|\mathcal G\|^2_2$ on the sparsity constraint $\mathcal{S}$. Using the vector of non-zero parameters $\phi=[\tilde{m}_1, \tilde{d}_1, \ldots, \tilde{m}_{n_c}, \tilde{d}_{n_c}]$ for the grid-following or alternatively $\phi=[\tilde{\alpha}_1, \tilde{\beta}_1, \ldots, \tilde{\alpha}_{n_c}, \tilde{\beta}_{n_c}]$ for the grid-forming implementation, the projected gradient is given by e.g., 
\[\text{proj}_{\phi}\left( \nabla_{K} \|\mathcal G\|^2_2\right) = \left(\tfrac{\partial}{\partial {\tilde{m}}_1} \|\mathcal G\|^2_2,\tfrac{\partial}{\partial {\tilde{d}}_1} \|\mathcal G\|^2_2, \ldots, \tfrac{\partial}{\partial {\tilde{d}}_{n_c}} \|\mathcal G\|^2_2\right).\]
Similar projections can be performed for the constraint set $\mathcal{C}$.

Because the $\mathcal H_2$ norm is infinite for unstable systems, both the system norm $\|\mathcal G\|^2_2$ as well as its gradient \eqref{eq:H2grad} are only well defined for a stable closed-loop system \eqref{eq:sysH2cl}. Thus, to optimize the control gain ${K}$, an initial guess for ${K}$ is required that stabilizes \eqref{eq:sysH2cl} and satisfies the constraints $\mathcal{S}$ and $\mathcal{C}$. Assuming that the system without VI devices is stable, it follows that $\tilde{m}_k=0$, $\tilde{d}_k=0$ suffices as an initial guess for the grid-following implementation. Moreover, the $\mathcal{H}_2$-norm is smooth and approaches infinity as the control gains $K$ approach the boundary of the set of stabilizing gains. In other words, any sequence of control gains $K$ with non-increasing cost is guaranteed to be stabilizing.

In the case that the projections onto $\mathcal{C}$ can be efficiently computed, the projected gradient method \cite{DPB:95} and gradient computation outlined above can be used to find a locally optimal solution to the optimization problem \eqref{eq:h2opt} even for systems of very large dimension. For instance this is the case when $\mathcal{C}$ encodes upper and lower bounds on $\tilde{m}_k$ and $\tilde{d}_k$. If the projection onto $\mathcal{C}$ cannot be computed efficiently, the above gradient computation can still be used to speed up the computation times of higher-order methods. 

\bibliographystyle{IEEEtran}
\bibliography{inertia}
\end{document}